\documentclass[11pt,twoside]{amsart}

\usepackage{amssymb}
\usepackage{amsxtra}
\usepackage{graphicx}
\usepackage{epsfig}
\usepackage{dcpic,pictexwd}

\theoremstyle{plain}

\newtheorem{thm}{Theorem}
\newtheorem{theorem}{Theorem}[section]
\newtheorem{prop}[theorem]{Proposition}

\newtheorem{defn}[theorem]{Definition}

\newtheorem{cor}{Corollary}

\newcommand{\al}{\alpha}
\newcommand{\be}{\beta}
\newcommand{\ga}{\gamma}

\newcommand{\Q}{\mathbb{Q}}
\newcommand{\R}{\mathbb{R}}
\newcommand{\T}{\mathbb{T}}

\newcommand{\Z}{\mathbb{Z}}

\newcommand{\Si}{\Sigma}

\newcommand{\ben}{\begin{enumerate}}
\newcommand{\een}{\end{enumerate}}

\newcommand{\SP}[1]{\mathfrak{#1}}
\newcommand{\SPS}[2]{\mathfrak{#1}_{#2}}
\newcommand{\ra}{\rightarrow}
\newcommand{\lra}{\longrightarrow}

\newcommand{\Oz}{P.\ Ozsv{\'a}th\,}
\newcommand{\Sz}{Z.\ Szab{\'o}\,}

\newcommand{\pic}[1]{\parbox{.6cm}{\psfig{figure=#1.eps,height=.6cm}}}

\begin{document}
\title{On knot Floer homology for some fibered knots}
\author[L. Roberts]{Lawrence Roberts}
\address{Department of Mathematics, Michigan State University,
East Lansing, MI 48824}
\email{lawrence@math.msu.edu}
\thanks {The author was supported in part by NSF grant DMS-0353717 (RTG)}
\maketitle

\section{Introduction}

\noindent Consider a link in $S^{3}$ for which one specified component, the axis, is an unknot.  We denote this link by $B \cup \mathbb{L}$ where $B$ is the axis of $\mathbb{L}$ and assume throughout that $\mathbb{L}$ intersects the spanning disc of $B$ in an odd number of points. For example, \\

\begin{center}
\includegraphics[scale=0.5]{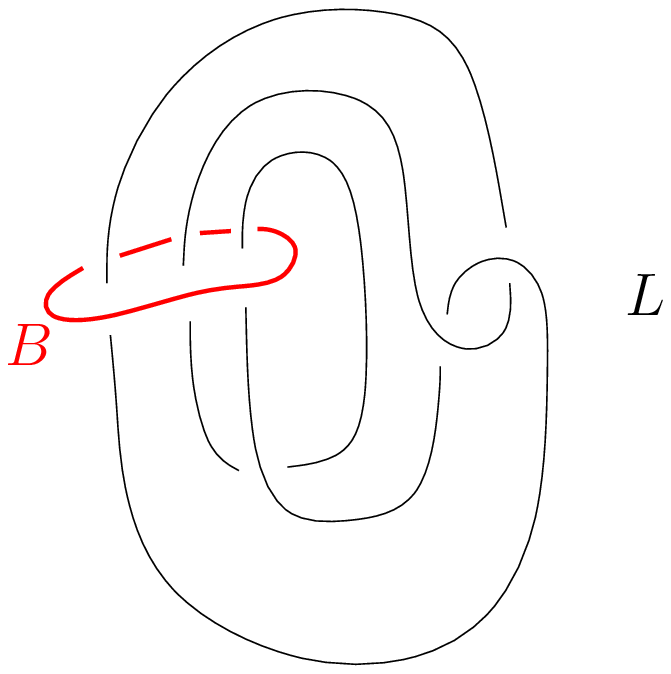}
\end{center}

\noindent Let $\Sigma(\mathbb{L})$ be the branched double cover of $S^{3}$ over $\mathbb{L}$, and let $\widetilde{B}$ be the pre-image of $B$ in $\Sigma(\mathbb{L})$. Then $\widetilde{B}$ is a null-homologous knot in $\Sigma(\mathbb{L})$ and we can try to compute 
$$
\widehat{HFK}(\Sigma(\mathbb{L}), \widetilde{B}, i) = \bigoplus_{\{\overline{\SP{s}}\, |\, \langle c_{1}(\underline{\SP{s}}), [F] \rangle\, =\, 2i\}}
\widehat{HFK}(\Sigma(\mathbb{L}), \widetilde{B}, \underline{\SP{s}}) 
$$
where $\underline{\SP{s}}$ is a relative $Spin^{c}$ structure for $\widetilde{B}$ and $[F]$ is the homology class of a pre-image of a 
spanning disc for $B$. The author began studying this situation in \cite{Robe} where a connection to Kohvanov homology is described. In this
paper, we wish to use the same approach to study the special case where $\mathbb{L}$ is a braid in the complement of $B$, and derive the complete
knot Floer homology for a myriad of fibered knots. \\
\ \\
\noindent We start by revisiting the main result of \cite{Robe} and proving it in a purely Heegaard-Floer manner, in order to use $\Z$-coefficients rather than $\Z/2\Z$-coefficients as in \cite{Robe}. For $\mathbb{L}$ alternating for the projection $A \times I \ra A$ this yields 

\begin{thm}
Let $\mathbb{L}$ be a non-split, alternating link in $A \times I$, with $\mathrm{det}(\mathbb{L}) \neq 0$, and which intersects the spanning disc for $B$ in an odd number of points. Then the $\Z/2\Z$-graded knot Floer homology $\oplus_{i\in\Z}\widehat{HFK}(\Si(\mathbb{L}), \widetilde{B}, F, i)$ is determined by a certain Turaev torsion, $\check{\tau}(\Si(\mathbb{L}) - K)$. Furthermore, for each $Spin^{c}$ structure, $\SP{s}$, on the $L$-space $\Si(\mathbb{L})$ we have $\tau(\widetilde{B}, \SP{s}) = 0$.
\end{thm}
\ \\
\noindent First, we review the definition of $\check{\tau}$ from \cite{Grig}. We will then use the above theorem to analyze the Heegaard-Floer homology of fibered three-manifolds whose monodromies can be represented as branched double covers of alternating braids. This is followed by several examples which should clarify the approach and can be read after the proof of the main theorem. In the final section, we concentrate on deriving results
in Heegaard-Floer homology similar to those of P. Seidel, I. Smith, and especially E. Eftekhary, \cite{Efte}, for the Floer cohomology of symplectomorphisms of surfaces. In particular, we will prove

\begin{prop}
Let $M_{\phi}$ be the fibered three manifold determined by a fiber $F^{g}$ and monodromy $D_{1}^{n_{1}} \cdots D_{2g}^{n_{2g}}$ where $n_{i} \geq 0$
and the Dehn twists occur along the linear chain of loops depicted in section 4. Let $\mathcal{S}$ be a collection of loops in $F$ consisting of $n_{i}$
parallel copies of the loop $\gamma_{i}$ in the linear chain. Then
$$
HF^{+}_{\Z/2\Z}(M_{\phi}, \SPS{s}{g-2}) \cong H^{\ast}(F\backslash \mathcal{S})
$$
as $H^{\ast}(F, \Z/2\Z)$-modules, where the action is by cup product on the right side of the isomorphism and by the $H_{1}$-action on the left. The Heegaard-Floer group is the direct sum of the homologies over all $Spin^{c}$ structures pairing with the fiber to give $2g-4$.
\end{prop} 

\noindent By duality, there is a corresponding theorem when $n_{i} \leq 0$ for all $i$.

\section{Background on Alexander polynomials}

\noindent Let $Y$ be a rational homology sphere; and let
$K \hookrightarrow Y$ be a null-homologous knot with spanning surface $F$. Then $H_{1}(Y - K,\,\Z) \cong H_{1}(Y,\,\Z) \oplus \Z$. Let $G = \pi_{1}(Y-K)$. By duality, $F$ defines a cohomology class in $\phi : G  \ra \Z$. We let $\widetilde{X}$ be the $\Z$ covering determined by this
cohomology class, and let $A_{\phi} = H_{1}(\widetilde{X}, \widetilde{p})$. The Alexander polynomial, $\Delta_{\phi}$, is defined to be the greatest common divisor of the elements of the first elementary ideal of $A_{\phi}$. \\
\ \\
\noindent For $Y$ an $L$-space, the absolute $\Z/2\Z$ grading in Heegaard-Floer homology assigns each $\widehat{HF}(Y,\SP{s}) \cong \mathbb{F}$
to the even grading. Moreover, this absolute grading corresponds on $\widehat{CF}$ to that given by the local intersection number at a generator, ${\bf x}$, between the two totally real tori, $\T_{\al}$ and $\T_{\be}$. It is chosen to ensure that $\chi(\widehat{HF}(Y)) = \big|H_{1}(Y,\,\Z)\big|$ 
(in fact, for any rational homology sphere). If we choose a Heegaard decomposition of $Y$ subordinate to $K$ and use the presentation of $G$ it
provides, we can recover the Alexander polynomial above by Fox calculus relative to the map $\phi$. Since the local intersection numbers
of the totally real tori correspond to the signs in the determinant employed in the Fox calculus, this will also be the Euler characteristic
of the knot Floer homology $\widehat{HFK}(Y,\,K,\,F)$, taken over all $Spin^{C}$ structures and using the $\Z/2\Z$-grading.\\
\ \\
\noindent In fact, we may choose a unique Alexander polynomial for $\phi$ by requiring that
\ben
\item[] $\bullet$ $\Delta_{\phi}(1) = \big| H_{1}(Y,\,\Z) \big|$
\item[] $\bullet$ $\Delta_{\phi}(t^{-1}) = \Delta_{\phi}(t)$
\een
The first statement is true of the Euler characteristic of $\widehat{HFK}(Y,\,K,\,F)$ because of the Euler characteristic properties of $\widehat{HF}(Y)$. The second statment is true due to the identity, for torsion $Spin^{c}$ structures, found in \cite{Knot}:
$$
\widehat{HFK}_{d}(Y,\,K,\,\underline{\SP{s}}) \cong \widehat{HFK}_{d - 2m}(Y,\,K,\,J\underline{\SP{s}})
$$
where $J$ is conjugation and $m = \frac{1}{2} \langle c_{1}(\underline{\SP{s}}), F \rangle$ (assuming $K = \partial F$ as an oriented knot). 
On $Y$ the conjugation of relative $Spin^{c}$ structures also conjugates the underlying structures on $Y$; the symmetry of the Alexander
polynomial exists because we sum over all $Spin^{c}$ structures. \\
\ \\
\noindent In \cite{Grig} a refined torsion, $\check{\tau}(Y - K) \in \Q(Spin^{c}(Y)(T))$, is constructed using Turaev's formalism.
It has the properties that 
$$
(T - 1) \check{\tau}(Y - K) = \sum_{\SP{s} \in Spin^{c}(Y)} p_{\SP{s}}(T) \cdot \SP{s}
$$
for 
$$
p_{\SP{s}}(T) = \sum_{i \in \Z} \chi\big(\widehat{HFK}(Y, \SP{s}; K, i)\big) T^{i}
$$
where the Euler characteristic is taken according to the absolute $\Z/2\Z$-grading. In particular, this torsion allows
us to distinguish the individual $Spin^{c}$ structures at the expense of a substantial increase in computational difficulty.
It is related to the previous Alexander polynomial by
$$
\sum_{\SP{s} \in Spin^{c}(Y)} p_{\SP{s}}(T) = \Delta_{\phi}(T)
$$  
To construct this element we need the map $\epsilon: \pi_{1}(Y - K) \ra H_{1}(Y;\,\Z)$ and a 
cell complex decomposition of $Y - K$. This produces a presentation for the fundamental group to which one applies
Fox's free differential calculus using the homomorphism $\phi \otimes \epsilon$. In our case, we can obtain the right
side of the equality above by considering the free differentials for generators other than the one from a meridian. 
This gives a square matrix and eliminates the pesky $(T-1)$ factors in the torsion computations. For more details, consult \cite{Grig}.

\section{Improving to $\Z$-coefficients}

\noindent Our first goal is to give a proof of the following result, which is a more specific version of the theorem from \cite{Robe}. $A$ is a round annulus in $\R^{2}$ to fix the embedding of $A \times I$ in $S^{3}$

\begin{theorem}\label{thm:alex}
Let $\mathbb{L}$ be a non-split alternating link in $A \times I$ intersecting the spanning disc for $B$ in an odd number of points. Then for each $k$ 
the $\Z/2\Z$-graded homology $\widehat{HFK}_{\Z}(\Si(\mathbb{L}), \widetilde{B}; \SP{s}, k)$ has rank determined
by the corresponding coefficient of $(T-1)\check{\tau}(\Si(\mathbb{L}) - \widetilde{B})$. The knot Floer spectral sequence collapses at the
$E^{2}$ page for any $\SP{s} \in Spin^{c}(\Si(\mathbb{L}))$, and $\tau(\widetilde{B}, \SP{s})  = 0$. Finally, if we filter $CFK^{\infty}(\Si(\mathbb{L}), \widetilde{B})$ using $[{\bf x}, i, j] \ra i + j$, the induced spectral sequence also collapses at the $E^{2}$
page.
\end{theorem} 

\noindent {\bf Proof:} as in \cite{Robe} this is proved by induction on the number of crossings in $\mathbb{L}$. In particular, either of
the resolutions of a crossing of $\mathbb{L}$ results in an alternating link with fewer crossings to which the result should apply. These
resolutions correspond to two terms in a surgery exact sequence whose third term is the desired fibered knot. The homology of the last is 
isomorphic to the mapping cone of the former arising from the sequence. The consequences in general of this perspective are the subject of 
\cite{Robe} following in the footsteps of \cite{Doub}. We carry out the proof in a series of steps. For more detail  on the approach
see \cite{Robe}. \\
\ \\
\noindent {\bf \#I:} The base case of the induction. Consider knots of the form

\begin{center}
\includegraphics[scale=0.5]{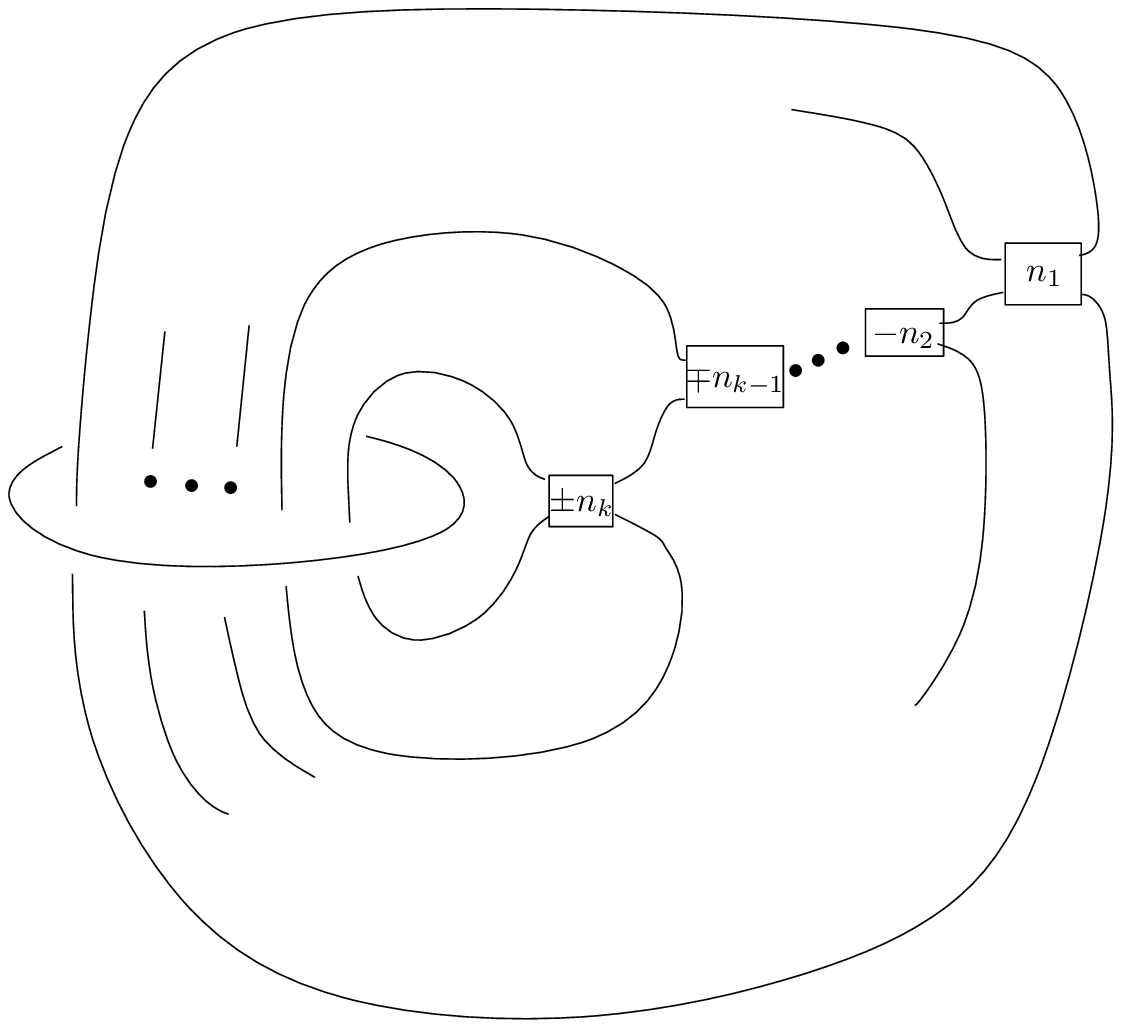}
\end{center}

Since $\mathbb{L}$ is an unknot, the branched double cover is $S^{3}$. $B \cup \mathbb{L}$ forms a link, which
can be described by flattening $B$ along a plane and pulling $\mathbb{L}$ so that the loops linking $B$ occur consecutively in
one direction. Thus we can untwist $\mathbb{L}$ at the expense of $B$, making $\mathbb{L}$ into an unknot. Choose an 
alternating projection for this link and consider the branched double cover. $\widetilde{B}$ will be an alternating knot
in $S^{3}$. The knot Floer homology has the following properties:
\ben
\item The grading of the knot Floer groups is determined by the filtration index according to $\mathcal{F} - \tau(K)$. In particular,
the $\Z/2\Z$ grading is just $\mathcal{F}$ modulo $2$. Thus the Alexander polynomial determines the knot Floer homology.
\item Consequently, the knot Floer spectral sequence collapses after the $E^{2}$ page. In fact, the $i + j$ spectral
sequence for $CFK^{\infty}$ collapses after the $E^{2}$ page.
\item $\tau(K) = -\frac{1}{2} \sigma(K)$ where $\sigma(K)$ is the signature of the knot
\een
\ \\
\noindent {\bf \#II:} We can calculate $\tau(\widetilde{B})$ for the twisted unknots. Put in the form alternating $B$ with $\mathbb{L}$ as the
axis (reverse the process above, interchanging the roles of the two components). Take the mirror if necessary, so
that the outermost region for the projection of $B$ will be colored black according to our coloring convention. Then $0 = \sigma(B)
= O(D) - 1 - n_{+}$ where $O(D)$ is the number of black regions, \cite{ELee}. A projection for $\widetilde{B}$ can be obtained by stacking
two copies of the tangle picture for $B$ determined by $\mathbb{L}$ and taking the closure. This has $2O(D) - 1$ black regions (since 
the outermost does not get doubled) and $2n_{+}$ positive crossings. Thus $\sigma(\widetilde{B}) = 2(O(D) - 1 - n_{+}) = 0$. So, our base cases all satisfy the required set of properties. \\
\ \\
\noindent {\bf \#III:}
We now follow Wehrli's algorithm, \cite{Wehr}, considered in the branched double cover. Recall that Wehrli's algortihm starts by enumerating the crossings of $\mathbb{L}$. One then proceeds through the crossings in order, looking at the two resolutions. If neither resolution disconnects the underlying four valent graph determined by $\mathbb{L}$, we resolve in both ways and then proceed to iterate the algorithm on the two resolved diagrams. To be specific we call the resolutiond the $0$ and the $1$ resolution where
$$
\pic{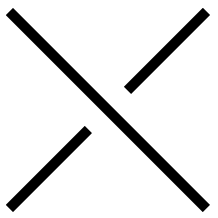}\stackrel{0}{\lra}\pic{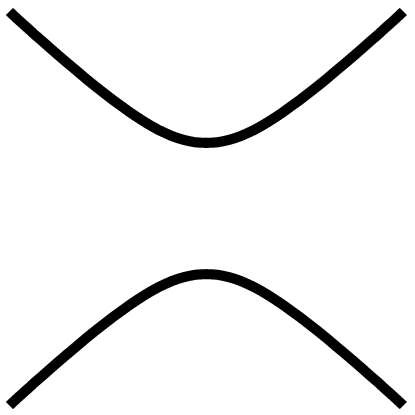} \hspace{.75in}
\pic{Lp}\stackrel{1}{\lra}\pic{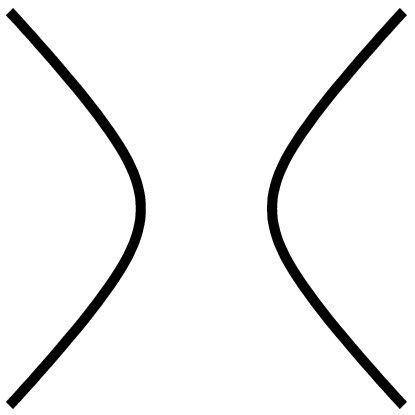}
$$

\noindent If either resolution disconnects the graph we move on to the next crossing. This produces a tree of resolutions whose leaves are used in \cite{Wehr} to give a smaller complex for calculating the Khovanov homology. Each leaf is associated to a spanning tree for the Tait graph of $\mathbb{L}$. In the double cover crossings correspond to $\pm 1$ fiber framed surgeries
on a specific loop in the fiber upstairs. One resolution contributes nothing to the monodromy. The other resoolution introduces a pair of critical points into the $S^{1}$-valued Morse function. We can depict this process downstairs by resolving $\mathbb{L}$. From the knot Floer homology
surgery sequence applied to each crossing according to the algortihm, we get a surgery spectral sequence, \cite{Robe}.  \\
\ \\
\noindent {\bf \#IV:}
To complete the argument we consider a single surgery curve and consider the homology of the $\pm 1$ surgery as
a mapping cone of the map between the other two. By considering the Euler characterstics in the associated long exact sequence, we have
the following relationship for the coefficients in the Alexander polynomial:
$$
a^{+1}_{j} = \pm a^{\infty}_{j} \pm a^{0}_{j}
$$
A similar argument for $-1$ surgery gives
$$
a^{-1}_{j} = \pm a^{\infty}_{j} \pm a^{0}_{j}
$$
where the signs would be determined by which maps shift the $\Z/2\Z$-gradings, and do not depend on $j$. However, we assume that
$$
\sum_{j} a^{*}_{j} = |H_{1}(Y_{\ast})|
$$
by our convention on Alexander polynomials. Since our manifolds are branched double covers over alternating links,
by the $L$-space arguments of \cite{Doub}, we have that $|H_{1}(Y_{\pm 1})| = |H_{1}(Y_{\infty})| + |H_{1}(Y_{0})|$ (where we have reverted to the framing conventions of \cite{Doub}). This fact reflects the association between spanning trees for the Tait graphs of alternating links and those of their resolutions. Since the signs are the same for all $j$, the only way both of these can be true is if the signs are positive. \\
\ \\
\noindent {\bf \#V:} We know $rk_{\pm 1, j} \leq rk_{\infty,j} + rk_{0, j}$. On the other hand, by
induction the terms on the right are $(-1)^{j} a^{\infty}_{j}$ and $(-1)^{j}a^{0}_{j}$, and $rk_{\pm 1, j} \geq |a^{\pm 1}_{j}|$
Using the identities above produces $rk_{\pm 1, j} = (-1)^{j}a^{\pm 1}_{j}$ which in turn implies that $a^{\pm 1}_{j}$ is negative for 
odd $j$ and positive for even $j$. In particular, all the homology in the $j^{th}$ level is contained in the same $\Z/2\Z$ grading. In the 
mapping cone construction, we must have that the chain map giving rise to the cone has $E^{1}(f)$ trivial
in order for the ranks to add. Therefore, the homology of $\pm 1$ surgery is the direct sum of the other two.\\
\ \\
\noindent {\bf \#VI:} Also by the $L$-space arguments in \cite{Doub} the $Spin^{c}$ structures on the new manifold partition into
two sets, one from each of the other two three manifolds. Since $\tau = 0$ for all the $Spin^{c}$ structures in these other manifolds, it
is also $0$ for all the $Spin^{c}$ structures here. Thus, the $\Z/2\Z$-gradings are determined as in the rank argument above. \\
\ \\
\noindent {\bf \#VII:} Finally, we note that the maps involved induce maps on the knot Floer spectral sequence. If the new knot Floer spectral
sequence has non-trivial higher differentials, the induced isomorphisms (from the $E^{1}$ terms being isomorphic) will force
higher differentials in one or the other of the resolved three manifolds. By induction this does not happen. Thus the 
spectral sequence will collapse at $E^{2}$. In fact, the chain maps also induce maps on the $i + j$ spectral sequences from
$CFK^{\infty}$ (which converged in a finite number of steps). Again, by induction this rules out the possibility of new
higher differentials.$\Diamond$\\ 
\ \\
\noindent The argument parallels that of \cite{Robe}, where there is more detail, but yields an algorithm we employ in the examples below. 
In fact, the above argument  holds for a broader class of links,  similar to the
quasi-alternating links of \cite{Doub}. This is the {\it smallest} subset of links in $A \times I$, denoted $\mathcal{Q}'$, with the property that
\ben
\item The alternating, twisted unknots, linking $B$ an odd number of times, are in $\mathcal{Q}'$.
\item 
 If $L \subset A \times I$ is a link admitting a connected projection to $A$, with a crossing such that
 \ben
 \item[] $\bullet$ The two resolutions of this crossing, $L_{0}$ and $L_{1}$, are in $\mathcal{Q}'$ and are
 connected in $A$, and
 \item[] $\bullet$ $\mathrm{det}(L_{i}) > 0$ for $i = 0,1$ and $\mathrm{det}(L) = \mathrm{det}(L_{0}) + \mathrm{det}(L_{1})$
 \een
 then $L$ is in $\mathcal{Q}'$
\een
Then alternating $L$ are in $\mathcal{Q}'$, and the elements of $\mathcal{Q}'$ when considered in $S^{3}$ are elements of $\mathcal{Q}$. For
this class of knots Wehrli's algorithm terminates at the base cases of our induction, where one or other resolution will disconnect the diagram.
 
\section{Alternating braids and alternating mapping classes}\label{sec:fibered}

\noindent We apply the preceding theory when $\mathbb{L}$ is a braid. If $\mathbb{L}$ has $b$ strands, then the branched double cover is fibered by genus $\frac{1}{2}(b - 1)$ punctured surfaces. To specify the monodromies we will consider, let $\gamma_{1}, \ldots, \gamma_{b-1}$ be the curves
depicted as:

\begin{center}
\includegraphics[scale=0.6]{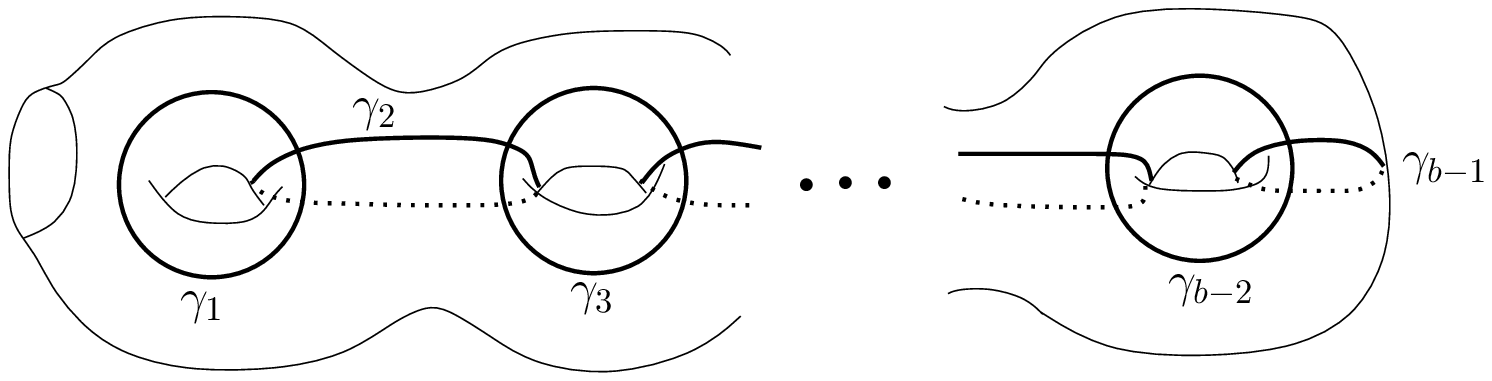}
\end{center}

\noindent Let $F$ be a surface of genus $g > 1$. Let $\mathcal{S} = \{ S_{1} , \ldots, S_{m}\}$ be simple closed loops in $F$. Assume that 
loops in $\mathcal{S}$ intersect once transversely or not at all. Define $G(\mathcal{S})$ to be the graph with vertices in one-to-one correspondence with the $S_{i}$ and edges in correspondence to the intersection points. This graph we shall call the {\it intersection graph} of $\mathcal{S}$. When
$G(\mathcal{S})$ is linear and the loops are non-separating, as above, we may use the theory in the previous sections to compute the knot Floer homology.

\begin{defn}
Let $\delta_{i}$ denote a positive Dehn twist around $\gamma_{i}$. An element $\phi \in \Gamma_{1}^{g}$ will be called {\em alternating} if it can be represented as a product $\delta_{i_{1}}^{\nu_{1}}\delta_{i_{2}}^{\nu_{2}} \cdots \delta_{i_{k}}^{\nu_{k}}$ where
\ben
\item If $i_{j} = i_{l}$ for some $j$ and $l$, then $\mathrm{sgn}(\nu_{j}) = \mathrm{sgn}(\nu_{l})$
\item If $i_{j} = i_{l} \pm 1$ for some $j$ and $l$ then $\mathrm{sgn}(\nu_{j}) = - \mathrm{sgn}(\nu_{l})$
\een
The element $\phi$ will be called {\em fully alternating} if there is such a representative for which 
$\{i_{1}, \ldots, i_{k}\} = \{1, \ldots, 2g\}$.
\end{defn} 

\noindent For an alternating representative of $\phi$, all the Dehn twists around a given circle will be performed in the same orientation. Any
$\phi \in \Gamma_{1}^{g}$ defines an open book decomposition of a three manifold and knot, $(Y,K)$, using a variation on the mapping torus
construction. We will mainly be concerned with the knot Floer homology of the binding for fully alternating mapping classes, but we can extend the results to some non-fully alternating mapping classes by taking connect sums of bindings and copies of $B(0,0)$, the knot in the connected sum of two $S^{1} \times S^{2}$'s obtained by performing $0$-surgery on two of the three components of the Borromean rings. \\
\ \\
\noindent Assume $b > 1$, so that $F$ has negative Euler chacteristic. The set $\gamma_{1}, \ldots, \gamma_{b-1}$ is a set of essential, simple closed curves which intersect efficiently and fill the surface (their complement is a boundary parallel annulus). For an alternating mapping class, the sets $\mathcal{G} = \{\gamma_{1}, \gamma_{3}, \ldots, \gamma_{b-2}\}$ and $\mathcal{D} = \{\gamma_{2}, \gamma_{4}, \ldots, \gamma_{b-1}\}$ satisfy the criteria for the main result of \cite{Penn}. $\phi$ may be reducible, but by \cite{Penn} each component map is either the identity or is pseudo-Anosov. If $\phi$ is fully alternating, then it is pseudo-Anosov.\\

\begin{cor}
Let $(Y, K)$ be a pair such that $Y$ is an open book with binding $K$, abstract page $F$, and fully alternating monodromy. Let $A$ be
the induced mapping on $H_{1}(F;\,\Z)$. If $\mathrm{det}(I - A) \neq 0$, then $Y$ is an $L$-space with $\mathrm{det}(I - A)$ $Spin^{C}$ 
structures. Furthermore, the $\Z/2\Z$-graded knot Floer homology is determined by $\check{\tau}(Y - K)$. For 
each $Spin^{c}$ structure on $Y$, $\tau(K, \SP{s}) = 0$.
\end{cor} 
\ \\
\noindent {\bf Proof:} $(Y,K)$ is the branched double cover of $(S^{3}, B)$ over an alternating braid $\sigma_{i_{1}}^{\nu_{1}}\sigma_{i_{2}}^{\nu_{2}} \cdots \sigma_{i_{k}}^{\nu_{k}}$. Since the monodromy is fully alternating, the branch locus is connected. By the second appendix the Alexander polynomial can be determined from $\mathrm{det}(I - tA)$, and $\mathrm{det}(I - A) = \pm \big|H_{1}(Y;\,\Z)\big|$. Thus $Y$
is a rational homology sphere. We now apply the results of \cite{Doub} to conclude that $Y$ is an $L$-space, and then apply theorem \ref{thm:alex} to compute the knot Floer homologies. $\Diamond$\\
\ \\
\noindent Let $\Delta = a_{0} + \sum_{i=1}^{n} a_{i}(T^{i} + T^{-i})$ be the Alexander polynomial for the fully alternating monodromy $\phi$. Define
the torsion coefficients for the binding, $\widetilde{B}$, by
$$
t_{s} =\sum_{j=1}^{\infty} j a_{|s| + j}
$$
then we can prove
\begin{prop}
Let $(Y,K)$ be an open book with fully alternating monodromy and pages of negative Euler characterisitic. Suppose further that $\mathrm{det}(I - A) \neq 0$. Let $Y_{K}$ be the fibered three manifold obtained by page framed surgery on $K$. Then, for all $s > 0$ we have a $\Z[U]$-module isomorphism
$$
\bigoplus_{\{\SP{s}\,:\,\langle c_{1}(\SP{s}), [\widehat{F}] \rangle = 2s \}} HF^{+}(Y_{K}, \SP{s}) \cong \Z^{b_{s}}
$$
where the the $\Z/2\Z$ grading of the right hand side is $(s\mathrm{\,mod\,}2)$, $[\widehat{F}]$ is the class in $H_{2}(Y_{K};\,\Z)$ 
for a capped page, and 
$$
b_{s} = (-1)^{s + 1} t_{s}(K)
$$ 
\end{prop} 
\ \\
\noindent {\bf Proof:} Since the knot Floer homology of the binding behaves like an alternating knot in $S^{3}$ we mimic the
proof of theorem 1.4 in \cite{Alte}. There are a few changes due to the lack of absolute grading information. 
First, note that for $p >> 0$, the isomorphism
$$
CF^{+}(Y_{p}(K), [s]) \lra C\{\mathrm{max}(i, j - s) \geq 0\} 
$$
from \cite{Knot} still applies in our setting, since we have a prescribed spanning surface. If we filter the left side by
$[{\bf x}, i] \ra i$ and the right side by $[{\bf y}, i, j] \ra i + j$, then both sides are filtered complexes. The chain isomorphism
above takes $[{\bf x}, i]$ to a sum of terms such as $[{\bf y}, i - n_{w}(\psi), i - n_{z}(\psi)]$ with $n_{w}(\psi) - n_{z}(\psi) = s$. 
Now $i - n_{w}(\psi) + i - n_{z}(\psi) = 2\,i + s - 2n_{w}(\psi) \leq 2i + s$. Thus the chain isomorphism is filtered and we have 
a spectral sequence morphism which converges to an isomorphism. For each $i + j$ value, the $E^{1}$ term on the right is a sum of
(shifted) knot Floer homology groups for the binding. From \ref{thm:alex} these have the property that
$$
\begin{array}{c}
gr_{\Z/2\Z}({\bf x}) = \mathcal{F}({\bf x}) \mathrm{\ mod}\,2\\
\ \\
gr_{\Z/2\Z}([{\bf x}, i, j]) = i + j  \mathrm{\ mod}\,2
\end{array}
$$
\ \\
\noindent The homology groups are constructed through the long exact sequence. At each stage of the long exact sequence there are chain maps
which we apply to $CFK^{\infty}$ to obtain spectral sequence morphisms for the $i + j$ filtration. If in $CFK^{\infty}$ there is a differential
which is not $(-1,0)$ or $(0,-1)$ then it must not induce a higher differential past the $E^{2}$ term. This occurs because the $E^{1}$ term
is isomorphic to a lower group in the resolution tree, and by induction these do not have higher differentials after the $E^{2}$ term (they
collapse at $E^{2}$ even in $HFK^{\infty}$). Thus, we need only consider up to the $E^{2}$ terms to calculate the homologies of these complexes
up to isomorphism. \\
\ \\
\noindent We can now proceed as in \cite{Alte}, replacing $C\{\mathrm{max}(i, j - s) \geq 0\}$, and every such complex, with $E^{1}\{\mathrm{max}(i, j - s) \geq 0\}$. Furthermore, we split these complexes up according to the $Spin^{c}$ structure on $Y$. We further decompose each of these $E^{1}$ pages
into subgroups $E^{1}_{k}$ for $k \in \Z$ by using those generators with absolute grading given by 
$$
gr_{\Q}([{\bf x}, i, j]) = i + j + d(\SP{s}) + k
$$
These are subcomplexes of $E^{1}(\SP{s})$ since only the $(-1,0)$ and $(0,-1)$ differentials will contribute. For each $E^{1}_{k}$ we apply the
argument from \cite{Alte}, noting that only for $k=0$ will we obtain a tower $\mathcal{T}^{+}$. Every other $k$ merely produces some finite group
in a specific grading ($s - 1 + k + d(\SP{s})$). The grading subscripts in the reduced homologies in \cite{Alte} will now only record the $\Z/2\Z$-grading. Now, however, we obtain for large positive surgeries on the binding, 
$$
E^{2}_{\SP{s}}\{\mathrm{max}(i, j - s) \geq 0\} \cong \mathcal{T}^{+}_{d(\SP{s})} \oplus \Z_{s - 1 \mathrm{\ mod}\,2}^{m_{(\SP{s},s)}}
$$
for each of the $Spin^{c}$ structures where $d(\SP{s})$ defines the $0$ in the $\Z/2\Z$-grading. So the reduced homology
in each of these homologies occurs in grading $s - 1$ modulo $2$. Since the spectral sequence for $CFK^{\infty}$ collapse at this
point, the tower in $E^{2}$ is the tower in $E^{\infty}$. The reduced homology all has the same $\Z/2\Z$-grading, and thus the 
spectral sequence collapses entirely. Using the long exact sequence,

$$
\cdots \lra HF^{+}(Y, \SP{s}) \stackrel{F_{1}}{\lra} HF^{+}(Y_{K}, \SP{s}_{s}) \stackrel{F_{2}}{\lra} HF^{+}(Y_{p}(K), \SP{s}_{s}) \stackrel{F_{3}}{\lra} \cdots
$$
\ \\
\noindent we know that $HF^{+}(Y, \SP{s}) \cong \mathcal{T}^{+}_{d(\SP{s})}$ is in the even absolute gradings, as is the tower in $HF^{+}(Y_{p}(K), \SP{s}_{s})
\cong E^{2}_{\SP{s}}\{\mathrm{max}(i, j - s) \geq 0\}$ with an absolute grading shift. The map $F_{3}$ is modelled on the surjection 
 $C\{\mathrm{max}(i, j - s)\} \ra C\{i \geq 0\}$. Since $\tau(\widetilde{B}) = 0$ and $s > 0$, it maps the tower surjectively onto the tower in $HF^{+}(Y, \SP{s})$. $F_{2}$ corresponds to a positive definite cobordism and therefore reverses the $\Z/2\Z$-gradings, implying that $F_{1}$ (a cobordism with a new $0$ framed homology class) and $F_{3}$ preserve them. Therefore, the $\Z_{s - 1\mathrm{\ mod}\,2}^{m_{(\SP{s},s)}}$ term gives rise
 to a $\Z^{m_{(\SP{s},s)}}$ in $\Z/2\Z$-grading given by $s$ modulo 2. \\
 \ \\
\noindent Adding these over all $Spin^{c}$ structures provides the identification with the torsion coefficients, by way of the result that
$-t(Y_{K},\SP{s}) = \chi(HF^{+}(Y_{K}, \SP{s}))$ since $b_{1}(Y_{K}) = 1$. Note that we use the fact that
the projection of $\tau(M)$ under projection to $\Z[H_{1}/\mathrm{Tors}]$ is still $\Delta_{K}/(t-1)^{2}$. In fact, the values
of $m_{(\SP{s},s)}$ should be given by the Turaev torsion directly. $\Diamond$
 
\section{Examples}

\begin{figure}
\begin{center}
\includegraphics[scale=0.6]{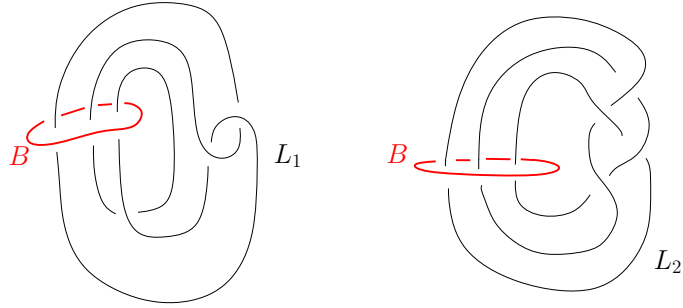}
\end{center}
\caption{The diagram for example 1 is on the left; that for example 2 is on the right.}\label{fig:examp1}
\end{figure}
 
\noindent {\bf Example 1:} See Figure \ref{fig:examp1} for the diagram. Here $\mathbb{L}$ is an unknot in $S^{3}$, so $\widetilde{B}$ is a knot
in $S^{3}$ as well. Untwisting and taking the branched double cover (or using symmetry between the two components) shows that $\widetilde{B}$ is 
the knot:

\begin{center}
\includegraphics[scale=0.6]{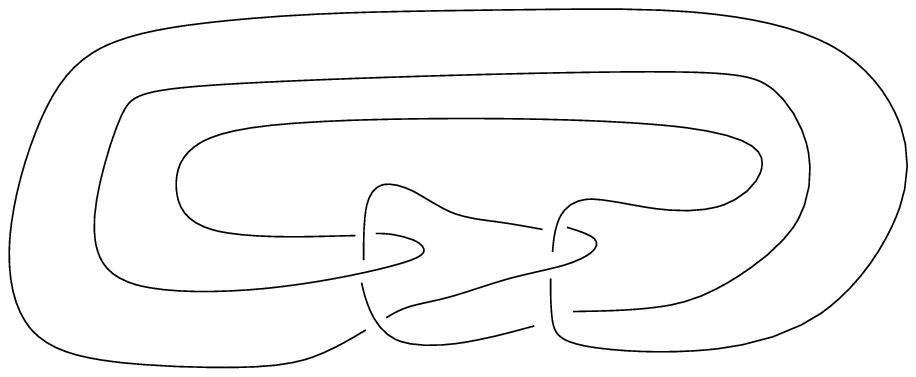}
\end{center}

\noindent This is the alternating knot, $6_{1}$, with signature equal to $0$. The main result in \cite{Alte} now verifies  the knot Floer conclusions of theorem \ref{thm:alex}. We note for later that the Alexander polynomial is $-2\,T^{-1} + 5 - 2\,T$. \\
\ \\
\noindent {\bf Example 2:} See Figure \ref{fig:examp1} for the diagram. Here $\mathbb{L}$ is the figure-8 knot, $4_{1}$, whose branched double cover is $L(5,2)$. In this arrangement, $\widetilde{B}$ is a genus $1$ fibered knot in $L(5,2)$. The possibilities for the homology of such a knot are strictly limited, since there is only a $\Z$ in filtration levels $\pm 1$. The real content of the theorem here is that $\tau(\widetilde{B}) = 0$, as this
implies that there is one $Spin^{c}$ structure where the knot Floer homology is that of $4_{1}$. We give a non genus $1$ example later.\\
\ \\
\noindent The monodromy for this knot is $\big(\gamma_{1}\gamma_{2}^{-1}\big)^{2}$ where
$\gamma_{i}$ is a positive Dehn twist around a standard symplectic basis element for $H_{1}(T^{2} - D^{2})$. The monodromy action on $H_{1}$ and the
Alexander polynomial associated to the $\Z$-covering from the fibering are computed to be 
$$
A = \left[\begin{array}{cc} 2 & 3 \\ 3 & 5 \\ \end{array} \right]  \hspace{.5in} \Longrightarrow \hspace{.5in} \mathrm{det}(I - tA) \doteq  \Delta_{\widetilde{B}}(t)  = - T^{-1} + 7 - T^{1}
$$
where we have symmetrized and normalized $\mathrm{det}(I - tA)$ according to our convention. In fact
we should use the more refined torsion, $\check{\tau}(Y - K)$, in our Euler characteristic computations, \cite{Grig}. This we now proceed to 
calculate.\\
\ \\
\noindent The fundamental group of $\Si(\mathbb{L}) - \widetilde{B}$ can be computed using the basis for $\pi_{1}(F)$ above. The action of $(D_{\gamma_{1}}D_{\gamma_{2}}^{-1})^{2}$ on the two elements generating this free group is
$$
\begin{array}{l}
\gamma_{1} \lra \gamma_{1}\gamma_{2}\gamma_{1}^{2}\gamma_{2}\gamma_{1}\gamma_{2}\gamma_{1} = R(\gamma_{1}) \\
\ \\
\gamma_{2} \lra \gamma_{1}\gamma_{2}\gamma_{1}^{2}\gamma_{2} = R(\gamma_{2})
\end{array}
$$
These provide the relations $\gamma_{1}^{-1}tR(\gamma_{1})t^{-1}$ and $\gamma_{2}^{-1}tR(\gamma_{2})t^{-1}$ for the fundamental
group. \\
\ \\
\noindent For the choice above, we obtain the map on homology $e_{1} \ra 5 e_{1} + 3 e_{2}$ and $e_{2} \ra 3 e_{1} + 2 e_{2}$. 
The quotient $\Z^{2}/(I - A)$ has a basis given by $(5,0)$ and $(3,1)$. The map to $H_{1}(\Si(\mathbb{L}))$ is thereby given as
$\gamma_{1} \ra e$ and $\gamma_{2} \ra e^{-3}$ (switching to exponents) with $e^{5} = 0$. We now apply
Fox calculus to the relations, and then map to $H_{1}(\Si(\mathbb{L}) - \widetilde{B})$ using the previous map and $t \ra T$. We  illustrate with
one calculation:
$$
(\rho \otimes \epsilon)(\partial_{\gamma_{2}}\,R_{1} ) = (\rho \otimes \epsilon)( \gamma_{1}^{-1}t\big( \gamma_{1} + \gamma_{1}\gamma_{2}\gamma_{1}^{2} + \gamma_{1}\gamma_{2}\gamma_{1}^{2}\gamma_{2}\gamma_{1} \big)) $$
$$
= e^{-1} t ( e + e\cdot e^{-3} \cdot e^{2} + e^{-2}) = t( 1 + e^{4} + e^{2})
$$
The overall matrix (removing the column for the derivatives related to $t$) is
$$
\left[ \begin{array}{cc}
	-e^{4} + T(2e^{4} + e + e^{2} + e^{3}) & T( 1+ e^{2} + e^{4}) \\
	T( e + e^{2} + e^{3}) & -e^{3} + T(e^{3} + e^{4}) \\
	\end{array} \right]
$$
which has determinant $-T - T\,e + (1 - 3T + T^{2})e^{2} - T\,e^{3} - T\,e^{4}$. Symmetrizing in $T$ produces the correct Euler characteristic
up to signs. This determines the ranks in each of the filtration indices for each of the $Spin^{c}$ structures.  \\
\ \\
\begin{figure}
\begin{center}
\includegraphics[scale=0.6]{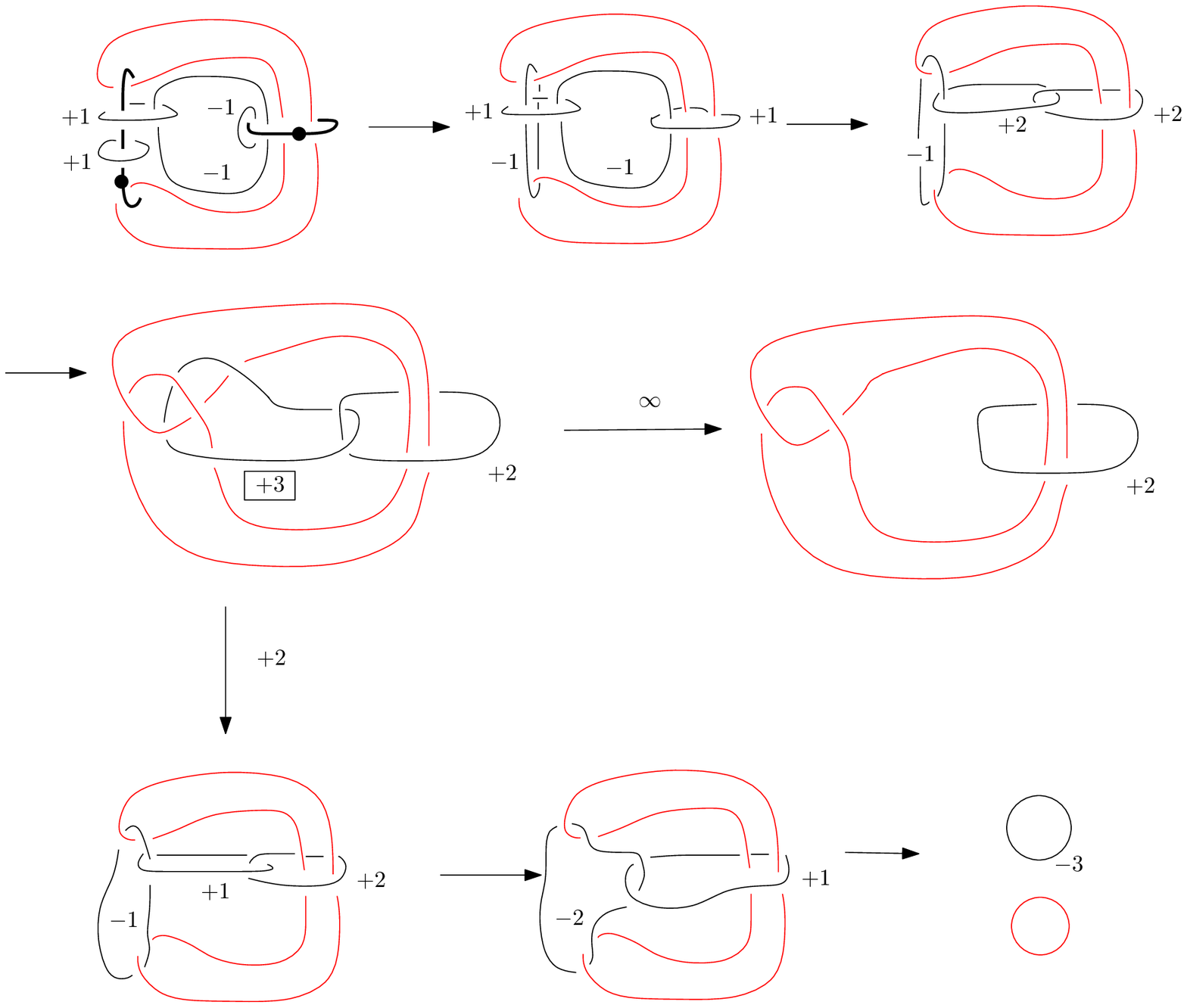}
\end{center}
\caption{ }\label{fig:L52}
\end{figure}

\noindent This fibered knot can also be seen by surgery on the diagram for $B(0,0)$ as in Figure \ref{fig:L52}. The box around $+3$ 
indicates that we will use it for the surgery sequence with framings $\infty, +2, +3$. For $\infty$
we obtain a knot in $L(2,1)$ identical with $B(-1,2)$. For $+2$ surgery we obtain the unknot in $L(-3,1)$. Since $L(5,2)$
has five $Spin^{c}$ structures, the sequence for the $\widehat{HFK}$ splits (all the homology from $L(-3,1)$ must map to that for
the fibered knot in $L(5,2)$). The knot $B(-2,1)$ has homology $\Z$  in filtration level $0$ for one $Spin^{c}$ structure, and
homology identical with $\widehat{HFK}(4_{1})$ in the other. This follows from a Borromean rings calculation or can be found
in \cite{Jabu}. Thus all the terms in filtrations other than $0$ should occur for a single $Spin^{c}$ structure (the one fixed under
conjugation). \\
\ \\
\noindent {\bf Example 3:} Consider the situation in section 5.2 of \cite{Jabu}: a genus 1 fibered knot with monodromy $D_{\gamma_{1}}^{n}D_{\gamma_{2}}^{m}$
with $m\cdot n< 0$. This is the branched double cover of the closure of the three stranded braid $\sigma_{1}^{n}\sigma_{2}^{m}$. Assume for now that $m < 0$. Then the action of the monodromy on $H_{1}(F)$ is given by 
   
$$
A = \left[\begin{array}{cc} 1 & n \\ 0 & 1 \\ \end{array} \right]\left[\begin{array}{cc} 1 & 0 \\ -m & 1 \\ \end{array} \right]  = 
\left[\begin{array}{cc} 1 - mn & -m \\ n & 1 \\ \end{array} \right] 
$$
\ \\
$$
\Longrightarrow \hspace{.5in}  \mathrm{det}(I - tA) \doteq \Delta_{\widetilde{B}}(t) = - T^{-1} + (2 - mn) - T^{1}
$$
By theorem \ref{thm:alex}, the Alexander polynomial determines the knot Floer homology groups. From Heegaard Floer homology, there must be an $\mathbb{F}$ in filtration level $0$ for each of the $\big|m\cdot n\big|$ $Spin^{c}$ structures on $\Sigma(\mathbb{L}) = L(m,1) \# L(n,1)$ (the closure of the braid is a connect sum of a $(2,m)$ torus link and a $(2,n)$ torus link). This implies that the knot Floer homology is given by
$$
\widehat{HFK}(\Si(\mathbb{L}), \widetilde{B}, j) \cong \left\{ \begin{array}{ccc} \mathbb{F} & \hspace{.25in} &  j = 1 \\ \mathbb{F}^{(2 - mn)} & \hspace{.25in} & j = 0 \\ \mathbb{F} & \hspace{.25in} & j = -1 \end{array} \right.
$$
with the $j=0$ level in the even grading. With the observation about absolute gradings in the proof of theorem \ref{thm:alex}, we recover Lemma 5.5 of \cite{Jabu} up to the decomposition into $Spin^{c}$ structures. \\
\ \\
\noindent Once again, we can compute $\check{\tau}(\Si(\mathbb{L}) - \widetilde{B})$. The map on $\pi_{1}(F)$ is
$$
\begin{array}{l}
\gamma_{1} \lra \big(\gamma_{1}^{n}\gamma_{2}\big)^{|m|}\gamma_{1} = R(\gamma_{1}) \\
\ \\
\gamma_{2} \lra \gamma_{1}^{n}\gamma_{2} = R(\gamma_{2})
\end{array}
$$
Following the procedure above produces
$$
(T-1)\check{\tau}(\Si(\mathbb{L}) - \widetilde{B}) = - T^{-1} + \left[ (1 + e_{1} \ldots + e_{1}^{n-1})(1 + e_{2} + \ldots + e_{2}^{|m|-1}) + 2 \right] - T
$$ 
where $e_{1}^{n} = e_{2}^{|m|} = 1$ and the map to $H_{1}(\Si(\mathbb{L})$ is given by $\gamma_{1} \ra e_{1}$ and $\gamma_{2} \ra e_{2}$. 
Thus, there is one $Spin^{c}$ structure where the knot Floer homology is that of $4_{1}$ and the rest are trivial. This is most of the result
in \cite{Jabu}. \\
\ \\
\noindent {\bf Example 4:} We outline an example for a pseudo-Anosov mapping class on a genus-two surface with boundary. Consider the five stranded 
braid $\sigma_{1}^{-2}\sigma_{3}^{-1}\sigma_{2}^{2}\sigma_{4}\sigma_{3}^{-1}$. This corresponds to the monodromy $D_{\gamma_{1}}^{-2}D_{\gamma_{3}}^{-1}D_{\gamma_{2}}^{2}D_{\gamma_{4}}D_{\gamma_{3}}^{-1}$. $\mathbb{L}$ then consists of a chain
of four unknots with linking numbers $-1$, $+1$, and $-1$ along the chain (using the standard braid orientation). We can compute that
the action of this mapping class on $\pi_{1}(F)$ is

\begin{figure}
\begin{center}
\includegraphics[scale=0.8]{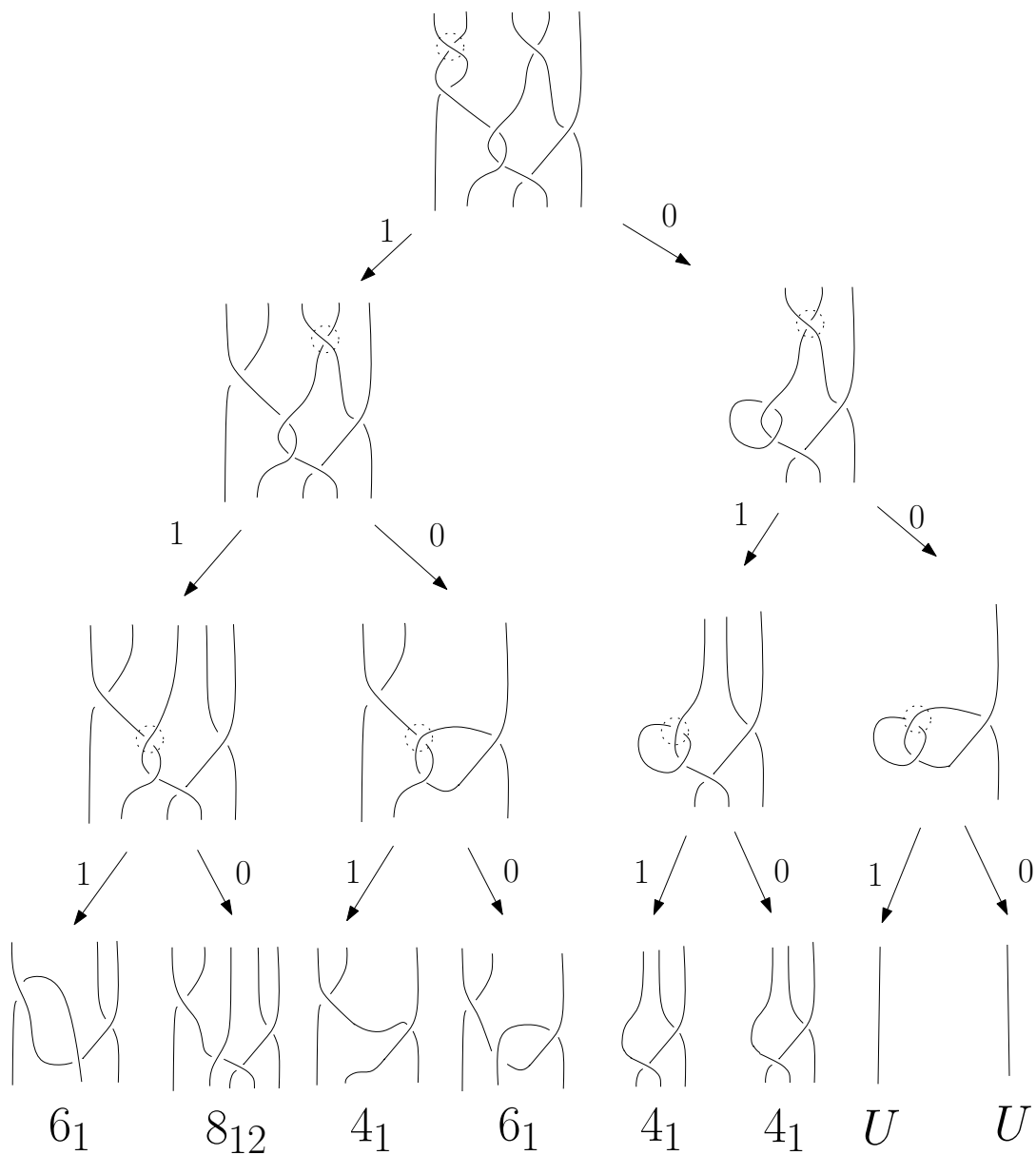}
\end{center}
\caption{The resolution tree for $\sigma_{1}^{-2}\sigma_{3}^{-1}\sigma_{2}^{2}\sigma_{4}\sigma_{3}^{-1}$ from Wehrli's algorithm. One should
think of these as standing in for their closures. We depict the results after resolving and simplifying. At each leaf we obtain a twisted unknot, and the label is the knot, $\widetilde{B}$, found in the branched double cover over this unknot (up to mirrors). Note that there is one $8_{12}$ label, two $6_{1}$ labels, three $4_{1}$ labels, and two unknots. This should be compared with the result in Example 4.}
\label{fig:restree}
\end{figure}

$$
\begin{array}{l}
\gamma_{1} \lra \gamma_{1}\gamma_{2}^{2} \\
\ \\
\gamma_{2} \lra \ga_{4}^{-1} \ga_{3}^{-1} \ga_{2} \ga_{1} \ga_{2} \ga_{4}^{-1} \ga_{3}^{-1} \ga_{2} \ga_{1} \ga_{2}^{2} \\
\ \\
\gamma_{3} \lra (\ga_{4}^{-1} \ga_{3}^{-1} \ga_{2} \ga_{1} \ga_{2})^{2} ( \ga_{3} \ga_{4} 
\ga_{2}^{-1} \ga_{1}^{-1} \ga_{2}^{-1})^{2}  \ga_{3} \ga_{4}^{2} 
(\ga_{3} \ga_{4} \ga_{2}^{-1} \ga_{1}^{-1} \ga_{2}^{-1})^{2}  \ga_{3} \ga_{4} \\
\ \\
\gamma_{4} \lra \ga_{4} \ga_{3} \ga_{4} \ga_{2}^{-1} \ga_{1}^{-1} \ga_{2}^{-1} \ga_{3} \ga_{4} \ga_{2}^{-1} \ga_{1}^{-1} \ga_{2}^{-1} \ga_{3} \ga_{4}\\
\ \\
\end{array}
$$
which yields the following map on $H_{1}(F)$: \\
$$
A = \left[ \begin{array}{cccc}
					1 & 2 & -2 & -2 \\
					2 & 5 & -4 & -4 \\
					0 & -2 & 4 & 3 \\
					0 & -2 & 5 & 4 \\
				\end{array} \right] \hspace{.25in} \Longrightarrow
\Delta_{\widetilde{B}}(T) \doteq T^{2} - 14 T + 34  - 14 T^{-1} + T^{-2} 
$$
\ \\
\noindent To distinguish the $Spin^{c}$ structures, we note that $\mathrm{det}(\mathbb{L}) = 8$ and $\Z^{4}/(A - I)$ $\cong \Z/2\,e_{1} \oplus \Z/2\,e_{2} \oplus \Z/2\,e_{3}$ with $e_{4} = e_{3}$ (where $e_{i} = [\gamma_{i}]$). Indeed, $\Si(\mathbb{L}) = L(2,1) \#^{2} L(-2,1)$. Calculating $(T-1)\check{\tau}(\Si(\mathbb{L}) - \widetilde{B})$ requires a 
great deal more effort, but ultimately yields:\\
$$
\begin{array}{l}
\big(T^{2} - 7 T + 13 - 7T^{-1} + 1\big) + \big(-2T + 5 - 2T^{-1}\big)\big(e_{1} + e_{3}\big) \\
\ \\
+ \big(-T + 3 - T^{-1}\big)\big(e_{2} + e_{1}e_{2}
+ e_{1}e_{3}\big) + \big(e_{2}e_{3} + e_{1}e_{2}e_{3}\big)
\end{array}
$$
\ \\
\noindent We make three observations: 1) setting $e_{i}$ to $1$ returns $\Delta_{\widetilde{B}}(T)$, 2) there are non-trivial phenomena in the knot Floer homology of more than one $Spin^{c}$ structure, and 3) the coefficient for each $Spin^{c}$ structure is the Alexander polynomial of a $\tau = 0$ alternating knot in $S^{3}$ (the first is that for $8_{12}$, the second for $6_{1}$, and the third for $4_{1}$). If one follows Wehrli's algorithm,
we obtain the tree of resolutions in Figure \ref{fig:restree}, which demonstrates how the knot Floer homology is built out of simpler pieces. \\
\ \\
\section{On a theorem of E. Eftekhary}

\begin{figure}
\begin{center}
\includegraphics[scale=0.6]{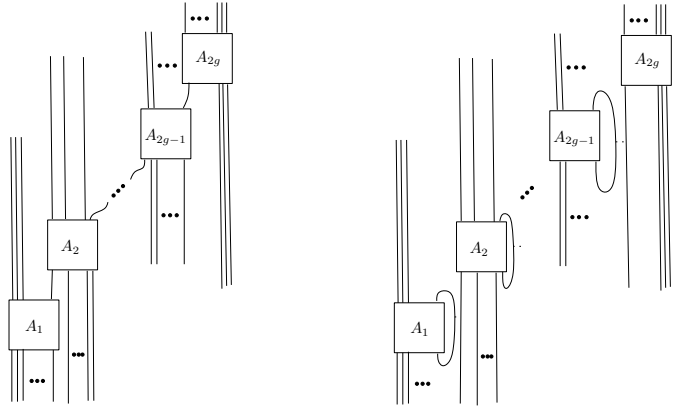}
\end{center}
\caption{Braids whose closures are connect sums. Each small loop in the right diagram is used to connect sum to the piece to its right. The result of all the connect sums is the picture on the left. If we assume that $A_{1}, \ldots, A_{2g}$ are connected alternating braids, each using an even number of strands, then the knot Floer homology of $\widetilde{B}$ for the left diagram is determined by the spanning trees of its Tait graph.}
\label{fig:connsum}
\end{figure}

\noindent In \cite{Efte} E. Eftekhary proves the following theorem for the Floer cohomology of symplectomorphisms:

\begin{thm}\cite{Efte}
Let $\mathcal{S}$ be a set of simple, closed, non-separating loops in a surface, $\widehat{F}$, each pair of which are either disjoint or intersect transversely in a single point and such that $G(\mathcal{S})$ is a forest. Let $\phi$ be the composition of a single positive Dehn twist along each loop in $\mathcal{S}$, taken in any order. Then
$$
HF^{\ast}_{symp}(\phi) = H^{\ast}(\widehat{F}, \mathcal{S})
$$
as $H^{\ast}(\widehat{F}, \Z/2\Z)$ modules where $H^{\ast}(\widehat{F})$ acts on the right side by the cup product and the left side by the quantum cup product.
\end{thm}
\noindent There is also a version for negative Dehn twists replacing $H^{\ast}(\widehat{F}, \mathcal{S})$ by $HF^{\ast}(\widehat{F}\backslash \mathcal{S})$ and a version for compositions of negative Dehn twists and positive Dehn twists as long as they occur on separated forests. Note that
we have an element of $\mathcal{S}$ for each Dehn twist; powers of the same Dehn twist should be construed as occuring along parallel copies, all of which are in $\mathcal{S}$, of a single curve. \\
\ \\
\noindent Due to the presumptive equivalence between various Floer homology theories, it has been suggested that a similar property should hold for
the Heegaard-Floer homology of the fibered three manifold in the above theorem. This approach questions whether the symplectic cohomology can be replaced by $HF^{+}(M_{\phi}, \SPS{s}{g-2})$, or some equivalent (using duality). Since $HF^{+}(M_{\phi}, \SPS{s}{g-1}) \cong \Z$ for every fibered three manifold, this is the next simplest case to try to compute. In \cite{Jabu}, the same statement using $HF^{+}$ is verified for certain
genus $1$ fibered three manifolds. Our purpose now is to extend their results to a certain case where $G(\mathcal{S})$ is a collection of linear chains
preserved by a hyperelliptic involution. We do this by first computing the knot Floer homology of the binding of an associated open book.\\
\ \\
\noindent Certain braid closures are connect sums of simpler pieces. These are depicted in Figure \ref{fig:connsum}, where each $A_{i}$ is a braid and the pieces are joined together in a staircase pattern. In order to have an odd number of strands intersecting the spanning disc for $B$, we ask that each $A_{i}$ intertwines an even number of strands. Furthermore, if these pieces
are alternating then the the number of spanning trees for each piece equals the number of $Spin^{c}$ structures on the double branched cover.
For the connect sums in Figure \ref{fig:connsum}, the number of trees is the product of the number of trees for each simpler piece, and 
the same is true for $Spin^{c}$ structures by the standard results for gluing along spheres. A similar conclusion to theorem \ref{thm:alex} then holds: the knot Floer homology of the branched double cover of the axis is a direct sum (over $Spin^{c}$ structures) of the knot Floer homologies of the double branched covers over the twisted unknots arising in the resolution tree. Furthermore,
we can determine the exact knot Floer homology by finding all the spanning tree unknots and examining the branched double covers (knots in $S^{3}$)
in those cases. This gives an effective algorithm for determing the knot Floer homology. However, now the value of $\tau(\widetilde{B})$ may
be different for the different pieces, since the base cases are no longer necessarily alternating. As one example consider the braids $\sigma_{1}^{n_{1}}\cdots \sigma_{2g}^{n_{2g}}$, where the exponents may be either positive or negative. When all the $n_{i} \geq 0$ we obtain a situation generalizing that of \cite{Efte} for mapping classes fixed by the hyperelliptic involution. The case when $g=1$ was addressed in \cite{Jabu}, and the technique in this paper recovers their results. For $g \geq 1$, tracing through the algorithm, using $0$ smoothings whenever possible, shows that there is one $Spin^{c}$ structure whose knot Floer homology will be identical with that of $T_{2,2g + 1}$. The others are more involved. \\

\ \\
\noindent For example, for $\sigma_{1}^{2}\sigma_{3}^{2}\sigma_{3}^{2}\sigma_{4}^{2}$, the Alexander polynomial is calculated from the monodromy action on $H_{1}(F)$ to be
$$
\Delta_{\widetilde{B}}(T) = T^{2} + 8\,T - 2 + 8\,T^{-1} + T^{-2}
$$
but the knot Floer homology, summed over all the $Spin^{c}$ structures of $\#^{4} L(2,1)$ is
$$
\widehat{HFK}(\widetilde{B}, j) = \left\{ \begin{array}{ll} 
\Z_{(0)} & j = 2 \\
\Z_{(-1)} \oplus \Z_{(0)}^{9} & j = 1 \\
\Z_{(-2)} \oplus \Z^{11}_{(-1)} \oplus \Z_{(0)}^{8} & j = 0 \\
\Z_{(-3)} \oplus \Z_{(-2)}^{9} & j = -1 \\
\Z_{(-4)} & j = -2 \\
\end{array} \right.
$$
where the gradings are $\Z$-relative gradings for the absolute grading in each $Spin^{c}$ structure. This computation follows by noting that
the resolution tree gives $8$ unknots, $5$ copies of $T_{2,3}$, 2 copies of $5_{2}$, and $1$ copy of $T_{2,5}$, as the branched double covers of the $16$ twisted unknots at the leaves. \\
\ \\
\noindent Now consider braids on $2g + 1$ strands of the form $\prod_{j} w_{j}$ where $w _{j} = \sigma_{i_{j}}\sigma_{i_{j} + 1} \ldots \sigma_{i_{j} + 2k_{j} - 1}$ and $i_{j+1} > i_{j} + 2k_{j}$. The last condition ensures that the braid words, $w_{j}$, include disjoint sets of generators. The monodromy of the open book in the branched double cover is then a series of negative Dehn twists along curves whose intersection graph forms a
forest of trees with no limbs. Furthermore at most one twist occurs along each circle. If we perform fibered framed $0$-surgery on the binding, the resulting fibered three manifold satisfies the conditions of the theorem in \cite{Efte} for the Floer cohomology of symplectomorphisms. To compute the knot Floer homology of $\widetilde{B}$ in this case, we will identify how
the knot in the branched double cover reflects the connect sum decomposition as just described. \\
\ \\
\noindent First, for each $\sigma_{j}$ not included in the product of the $w_{j}$'s there is a core circle, $c_{j}$, in the annulus, $A$, which does not intersect the diagram for $\mathbb{L}$, and which is the intersection of the plane with a sphere in $S^{3}$. This sphere can be chosen to intersect $B$ in two points. If we cut along all such spheres we obtain two $B^{3}$ pieces and a bunch of $S^{2} \times I$ pieces. Each piece contains some portion of $\mathbb{L}$ and either an arc (for the $B^{3}$'s) or two arcs (for the $S^{2} \times I$'s) of $B$. These pieces can be ordered from left to right in accordance with the ordering on the braid. We now fill every $S^{2}$ with a copy of $B^{3}$ and complete each set of arcs from $B$ with unknotted arcs in the new $B^{3}$ components to obtain a knot.  This realizes the pair $(B, \mathbb{L})$ as a pair connect sum along the axes of $(S^{3}, \mathbb{L}_{j}, B_{j})$ where $j$ indexes the pieces as in the braid. Furthermore, each $\mathbb{L}_{j}$ has an odd number of strands in it, due to the form of the braid word (we can alter this to included even number of strands and reach the same conclusion, but this requires more work). The branched double cover of $(S^{3}, \mathbb{L}_{j}, B_{j})$ over $\mathbb{L}_{j}$ is a copy of $(S^{3}, T_{2, 2k_{j} + 1})$ where a single strand of $\mathbb{L}$ gives rise to an unknot. \\
\ \\
\noindent We recover the branched double cover from the pieces in the following way. On $B_{i}$ there are one or two arcs which were glued in during the decomposition process. One arc for the pieces from either end of the braid; two for any piece from the interior. These lift to two or four arcs in the double cover (in the case of four, the lifts alternate around the knot $\widetilde{B}_{i}$). To construct the branched double cover, remove small ball neighborhoods around these and glue to the corresponding lift in the branched double cover of the pieces from the left and right in the braid diagram. One (or two) of these arcs glue just as connect sums of pairs (double cover, knot)to the piece on the left and/or right. To glue the other two arcs, add a four dimensional one handle and connect sum across the one handle. The boundary is the knot in the branched double cover. Note that this means that the knot sits in a connect sum of $S^{1} \times S^{2}$'s. \\
\ \\
\noindent We have made our choices so that there will be an even number of one handles added, say $2m$. From \cite{Knot}, in this setting
we can peel off $m$ copies of $B(0,0) \subset \#^{2} S^{1} \times S^{2}$ as connect summands, leaving a connect sum of the underlying knots $T_{2, 2k_{j} + 1}$ in $S^3$. This allows us to calculate the knot Floer homology of $\widetilde{B}$: $\widehat{HFK}(\widetilde{B}, g) \cong \Z_{m}$ while the computation of $\widehat{HFK}(\widetilde{B}, g - 1)$ is slightly more complicated. We have that $2m + \sum 2k_{j} = 2g$. We use the genus -1 level of each $B(0,0)$ or $T_{2, 2k_{j} + 1}$ in turn, tensored with the genus levels of all the others. Using the lower level in a $B(0,0)$ summand produces $$\Z_{0}^{2} \otimes \Z_{1} \otimes \cdots \otimes \Z_{1} \otimes \Z_{0} \otimes \cdots \otimes \Z_{0} = \Z_{m - 1}^{2}$$
There are m of these, for a total contribution $\Z_{m-1}^{2m}$. Using the lower filtration from a $T_{2, 2k_{j} + 1}$ summand produces 
$$
\Z_{1} \otimes \cdots \Z_{1} \otimes \Z_{-1} \otimes \Z_{0} \otimes \cdots \otimes \Z_{0}$$
contributing $\Z_{m-1}$ a total of $s$ times, where there are $s$ torus knots. Thus $\widehat{HFK}(\widetilde{B}, g-1) = \Z_{m-1}^{2m + s}$. \\
\ \\
\noindent We can now compute the Heegaard-Floer homology of $0$-surgery on $\widetilde{B}$. We use a theorem of \Oz and \Sz,
an account of which is given in \cite{Jabu}. Key to this approach is that the knot is in a three manifold without reduced homology (and only one $Spin^{c}$ structure with any non-trivial homology). This theorem states that $HF^{+}(Y_{\widetilde{B}}(0), \SPS{s}{g-2})$ is then isomorphic to the portion of the $\Z\oplus \Z$-graded complex $CFK^{\infty}$ described by $C\{i < 0, j \geq g-2\}$. For us this complex is isomorphic to

$$
\begindc{\commdiag}[4]
\obj(10,15)[T]{$\Z_{m - 2}$}
\obj(0,5)[t]{$\Z_{m - 4}$}
\obj(20,5)[B]{$ \oplus\ \Z_{m-3}^{2m}$}
\obj(10,10)[T1]{$\Z_{m-3}$}
\obj(10,6)[D]{$\vdots$}
\obj(10,0)[Ts]{$\Z_{m-3}$}
\obj(12,10)[E1]{$\ $}
\obj(12,0)[E2]{$\ $}
\mor{T1}{t}{$\ $}
\mor{Ts}{t}{$\ $}
\mor{E1}{E2}{$s$}[1,2]
\enddc
$$
where the arrows to the left are all isomorphisms. These arrows come from the complexes for the torus knots, specifically from the surjective differential onto the $-g(T_{2,2k_{j}+ 1})$ indexed summands. Thus $HF^{+}(Y_{K}, \SPS{s}{g-2}) \cong \Z_{m-2} \oplus \Z_{m-3}^{2m + s -1}$ where the
gradings should now be taken as relative gradings. Let $\widehat{F}$ be the genus $g$ fiber for $Y_{K}$. We have that $\widehat{F}\backslash C$ is a genus $g - \sum k_{j} = m$ surface with $s$ boundary components. Thus, $H^{\ast}(\widehat{F}\backslash C) \cong \Z \oplus \Z^{2m + s - 1}$ where the first summand is $H^{0}(\widehat{F}\backslash C)$ and the second is $H^{1}(\widehat{F}\backslash C)$.  Finally, the action of $H_{1}(Y_{K})/\mathrm{Tors}$, excluding the $\Z$ introduced during the surgery, is identical to that from $\#^{m} B(0,0)$, which corresponds (over $\Z/2\Z$) to the standard cup product on $\widehat{F}$. This verifies a Heegaard-Floer analog of E. Eftekhary's theorem for this set of braids. \\
\ \\
\noindent The theorem in \cite{Efte} is slightly more general, even in our setting, allowing words of the form $\sigma_{i_{j}}^{n_{j}}\sigma_{i_{j} + 1} \ldots \sigma_{i_{j} + 2k_{j} - 2}\sigma_{i_{j} + 2k_{j} - 1}^{m_{j}}$ with $n_{j}, m_{j} > 0$. It is a consequence of one of the lemmas in \cite{Efte}, or by inspection, that any reordering of this product can be brought into the form above by isotopies of the braid closure. Simply start at the left hand side, just above the first crossing in the second column and isotope all the crossings in the first column around the closed braid until they all lie under the crossing in the second column. Now take everything in the first and second column and isotope around until they all lie below the crossing in the third column. If we keep doing this we obtain a braid word of the kind at the beginning of this paragraph. In fact, this is all that can happen for a linear chain in his theorem . If the exponents are bigger than one in the interior of a word then $G(\mathcal{S})$ is no longer a forest.  \\
\ \\
\noindent Using our techniques we can generalize further, allowing words
where all the exponents can be arbitrary positive numbers: $$\prod_{j} \sigma_{i_{j}}^{n_{i_{j}}}\sigma_{i_{j} + 1}^{n_{i_{j} + 1}} \ldots \sigma_{i_{j} + 2k_{j} - 2}^{n_{i_{j} + 2 k_{j} - 2}}\sigma_{i_{j} + 2k_{j} - 1}^{n_{i_{j} + 2 k_{j} - 1}}$$
We now proceed to analyze this case.\\
\ \\
\noindent The closure of each word is subject to the analysis of the ladder braids mentioned at the beginning of this section. We need only determine the contributions of each $Spin^{c}$ structure on the branched double cover to $\widehat{HFK}(\widetilde{B}, k_{j})$ and $\widehat{HFK}(\widetilde{B}, k_{j} -1)$. As before, if we resolve each crossing using a $0$ resolution, when available, we obtain $T_{2,2k_{j}+1}$
in the branched double cover. This contributes a $\Z_{0}$ to $\widehat{HFK}(\widetilde{B}, k_{j})$ and a $\Z_{-1}$ to $\widehat{HFK}(\widetilde{B}, k_{j}-1)$. \\
\ \\
\begin{figure}
\begin{center}
\includegraphics[scale=0.6]{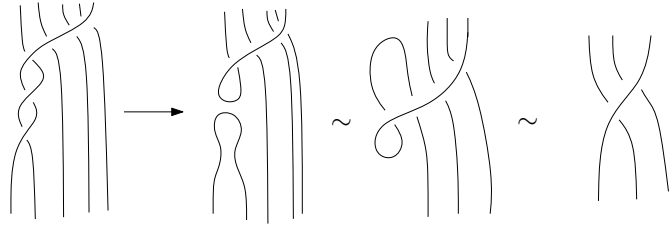}
\caption{If a one resolution is employed in either of the extreme twist regions we can isotope the closure of the fully resolved diagrams so
that there are two fewer strands.} \label{fig:oneoneres}
\end{center}
\end{figure}

\begin{figure}
\begin{center}
\includegraphics[scale=0.6]{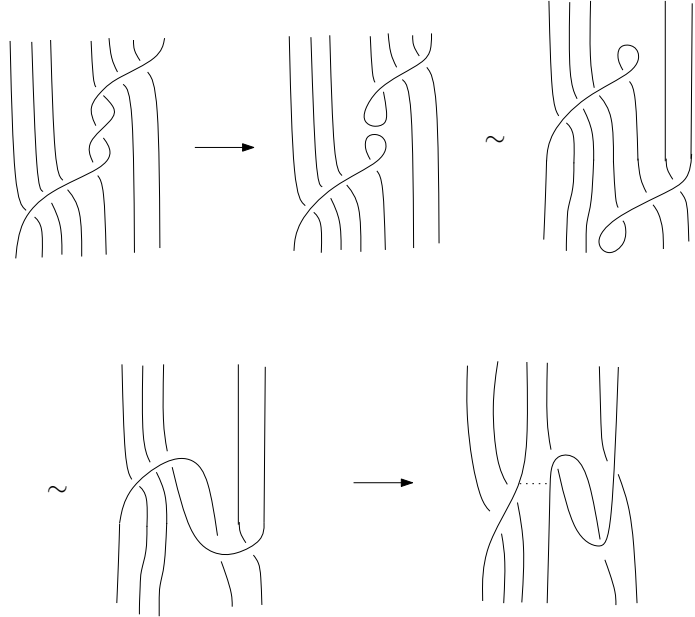}
\caption{If the single one resolution occurs in the middle of the braid, the closure is isotopic to the closure of a tangle with a clasp. The final diagram shows the pieces used in the plumbing decomposition of $\widetilde{B}$ in the branched double cover. The central piece alone will yield a copy of $5_{2}$ for $\widetilde{B}$ in the double branched cover, and this has a $\Z^{2}_{0}$ in its top knot Floer homology filtration level. Note that this is a positive braid, so we know that it is $5_{2}$ and not its mirror.} \label{fig:oneoneres2}
\end{center}
\end{figure}

\noindent Allowing a single $1$ resolution, and requiring all the others to be $0$ resolutions, contributes a copy of $T_{2, 2k_{j} -1}$ if the resolution occurs in either of the extreme twist regions, see Figure \ref{fig:oneoneres}. There are $n_{i_{j}} - 1$ and $n_{i_{j} + 2k_{j} -1} - 1$ such resolutions, respectively,  since we must still have a connected diagram when we use the $0$ resolution for the same crossing. If the $1$ resolution occurs in one of the interior twist regions we obtain a diagram as in Figure \ref{fig:oneoneres2} which contributes $\Z_{0}^{2}$ to $\widehat{HFK}(\widetilde{B},k_{j}-1)$. There are $n_{i_{j} +l} - 1$ of these occurring in the $i_{j} + l$ column. To see the contibution, frist intercahnge the roles of $\widetilde{B}$ and $\mathbb{L}$ and note that we can group an even number of crossings from each end of the braid so as to leave precisely three crossings, two for the clasp and one more to one or other side of the clasp. In the double cover, this grouping corresponds to viewing $\widetilde{B}$ as the plumbing of up to three objects, two torus knots $T_{2, 2k + 1}$ and $5_{2}$. In addition, the plumbing occurs along
genus minimizing spanning surfaces for each of the knots. By the Y. Ni's theorem, \cite{YiNi}, concerning Murasugi sums of knots, this gives the contribution as $\Z^{2}_{0}$, where $\widetilde{B}$ being a positive knot implies $\tau(\widetilde{B}) = g_{3}(\widetilde{B})$, \cite{Livi}, and thus confirms the grading.   \\

\begin{figure}
\begin{center}
\includegraphics[scale=0.6]{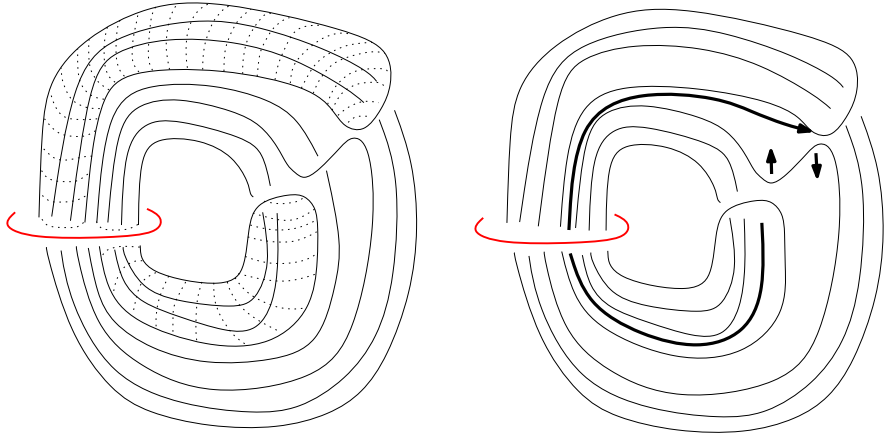}
\end{center}
\caption{The diagram on the left depicts the situation for two non-consecutive $1$ resolutions. The dotted lines trace the discs which lift to 
a compression disc for the fiber in the branched cover. The crossing assumptions ensure that the discs are above the diagram. The arcs in the spanning
disc for $B$ won't intersect since there is at least one strand between the two resolutions. The right diagram depicts the situation when the resolutions are consecutive. Now the compression discs will intersect in a point. Following the arrows will cancel the critical points in the diagram
and leave a copy of $\sigma_{1} \ldots \sigma_{2k}$.}
\label{fig:noncons}
\end{figure}

\ \\
\noindent When we allow two $1$ resolutions not occuring in consecutive columns, or three or more $1$ resolutions, we can see that there is no contribution to the $k_{j} - 1$ filtration level. Each $1$ resolution must occur in distinct columns, and there must be one such resolution that is the rightmost (and highest) and one which is the leftmost (and lowest). These are non-consecutive by assumption. As in Figure \ref{fig:noncons}, if we follow the strands opening down from the leftmost, and the strands opening up from the rightmost both will intersect the spanning disc for $B$ in cancelling pairs of points. Moreover, since these are not consecutive, the structure of the braid ensures that they will be disjoint. By following a small arc between these strands starting at the critical point, we obtain two arcs, necessarily disjoint, and two discs swept out by these arcs. These discs do not intersect $\mathbb{L}$ in their interiors and lift to compression discs for $F$ in the double cover. Hence the branched double cover for the resolved diagram has minimal genus for its binding smaller than $k_{j}-1$. \\
\ \\
\noindent This leaves only the case of two $1$ resolutions occuring in consecutive columns, see Figure \ref{fig:noncons} again. The two critical points from one resolution cancel with the two critical points from the other resolution. Hence we get a copy of $T_{2, 2k_{j} - 1}$ which contributes a $\Z_{0}$ to the $k_{j} - 1$ level. There are $(n_{i_{j} + l} - 1)(n_{i_{j} + l + 1} - 1)$ such contirbutions from the $i_{j} + l$ and $i_{j} + l +  1$ columns. Thus, adding up all these contributions yields that $\widehat{HFK}(\widetilde{B}, k_{j}-1)$ is congruent to $\Z_{-1} \oplus \Z_{0}^{T_{j}}$ where $T_{j} = \sum_{l} (n_{i_{j} + l} n_{i_{j} + l + 1} - 1 )$ since the latter is equal to 
$$
(n_{i_1} - 1) + (n_{i_{1}} - 1)(n_{i_{1}+1} - 1) + 2(n_{i_{1}+ 1} - 1) + (n_{i_{1} + 1} - 1)(n_{i_{1}+2} - 1) + \ldots + (n_{i_{1}+2k_{j} -1} - 1)
$$
\ \\
\noindent We now return to the product of all these braid words. As before $\widehat{HFK}(\widetilde{B}, g) = \Z_{m}$ but $$\widehat{HFK}(\widetilde{B}, g-1) \cong \Z_{m-1}^{2m} \oplus \Z_{m-1}^{s} \oplus \Z_{m}^{\sum_{j} T_{j}}$$ As before, the differential in the relevant subset of $CFK^{\infty}$ carries each
  generator of $\Z^{s}$ to $\Z$, coming from the spectral sequence for $T_{2, 2k_{j} + 1}$. So, we obtain $$HF^{+}(Y_{K}, \SPS{s}{g-2}) \cong \Z_{m}^{\sum \,T_{j} + 1} \oplus \Z_{m-1}^{2m + s - 1}$$ with action of $H_{1}$ entirely contained in the $\Z_{m}^{1} \oplus \Z_{m-1}^{2m}$
  portion of the complex (arising from the $B(0,0)$ summands). $\widehat{F}\backslash C$ consists of the surface in the simpler case along with
  numerous squares. It is straightforward to verify that there are $\sum T_{j}$ such  ``squares'' which each have rank one $H^{0}$-group. Again we also obtain the correspondence of the action of $H_{1}$ with the cohomology ring. \\
\ \\
\noindent Unfortunately, not all positive braids possess the property in Eftekhary's theorem. The braid $\sigma_{1}\sigma_{2}\sigma_{3}\sigma_{4}\sigma_{5}\sigma_{4}\sigma_{6}$ describes an open book upon which $0$ surgery has a rank $2$ $HF^{+}$-group for the relevant $Spin^{c}$ structure. This follows from the following calculation for closed braids
$$
\begin{array}{l}
C(\sigma_{1}\sigma_{2}\sigma_{3}\sigma_{4}\sigma_{5}\sigma_{4}\sigma_{6}) = C(\sigma_{1}\sigma_{2}\sigma_{3}\sigma_{5}\sigma_{4}\sigma_{5}\sigma_{6}) = \\
C(\sigma_{5}\sigma_{1}\sigma_{2}\sigma_{3}\sigma_{4}\sigma_{5}\sigma_{6}) = C(\sigma_{1}\sigma_{2}\sigma_{3}\sigma_{4}\sigma_{5}\sigma_{6}\sigma_{5}) =\\
C(\sigma_{1}\sigma_{2}\sigma_{3}\sigma_{4}\sigma_{6}\sigma_{5}\sigma_{6}) = C(\sigma_{6}\sigma_{1}\sigma_{2}\sigma_{3}\sigma_{4}\sigma_{5}\sigma_{6}) =
\\
C(\sigma_{1}\sigma_{2}\sigma_{3}\sigma_{4}\sigma_{5}\sigma_{6}^{2})
\end{array}
$$
The results above apply to the last braid closure, and yield that $HF^{+}(Y, \SPS{s}{1}) \cong 2$. However, $H^{\ast}(\widehat{F}\backslash C) \cong \Z^{3}$ for the first braid. The singular homology does not transform appropriately under the braid equivalences for the result to hold. From the proof above, we can identify the difficulty: the single allowable $1$ resolution for the original braid, although it does not occur at the ends of the braid, still produces a torus knot. The additional $\sigma_{4}$ allows the two critical points introduced during the resolution to cancel with each other instead of forming a clasp. This suggests first assuming some kind of normal form for our braids before proceeding further.

\end{document}